
\input amstex
\input epsf
\documentstyle{amsppt}
\nologo
\TagsOnRight
\def\epsfsize#1#2{\hsize}
\magnification 1200
\voffset -1cm

\def\fig#1{\bigskip\centerline{\epsffile{FIG2-#1.EPS}}\nopagebreak\medskip
\nopagebreak\centerline{\bf Fig\. #1}\bigskip}
\def\cal#1{\expandafter\def\csname#1\endcsname{\Cal{#1}}} \cal W \cal M
\def\ro#1{\expandafter\def\csname#1\endcsname{\operatorname{#1}}} \cal K
\ro{rk} \ro{Fr} \ro{Int} \ro{lk} \ro{im} \ro{ker} \ro{coim} \ro{Tors}
\let\emb\hookrightarrow \let\x\times \let\but\setminus \let\eps\varepsilon
\def\section#1{\bigskip\leftline{\bf #1}\nopagebreak\medskip}
\def\Cl#1{{\overline{#1}}} \def\N{\Bbb N} \def\Z{\Bbb Z} \def\R{\Bbb R}

\topmatter 

\title $n$-Quasi-isotopy: II. Comparison \endtitle
\author Sergey A. Melikhov and Du\v{s}an Repov\v{s}\endauthor

\abstract
Geometric aspects of the filtration on classical links by $k$-quasi-isotopy
are discussed, including the effect of Whitehead doubling and relations with
Smythe's $n$-splitting and Kobayashi's $k$-contractibility.
One observation is: $\omega$-quasi-isotopy is equivalent to PL isotopy for
links in a homotopy $3$-sphere (resp\. contractible open $3$-manifold) $M$
if and only if $M$ is homeomorphic to $S^3$ (resp\. $\R^3$).
As a byproduct of the proof of the ``if'' part, we obtain that every compact
subset of an acyclic open set in a compact orientable $3$-manifold $M$ is
contained in a PL homology $3$-ball in ~$M$.

We show that $k$-quasi-isotopy implies $(k+1)$-cobordism of Cochran and Orr.
If $z^{m-1}(c_0+c_1z^2+\dots+c_nz^{2n})$ denotes the Conway polynomial of
an $m$-component link, it follows that the residue class of $c_k$ modulo
$\gcd(c_0,\dots,c_{k-1})$ is invariant under $k$-quasi-isotopy.
Another corollary is that each Cochran's derived invariant $\beta^k$ is also
invariant under $k$-quasi-isotopy, and therefore assumes the same value on
all PL links, sufficiently $C^0$-close to a given topological link.
This overcomes an algebraic obstacle encountered by Kojima and Yamasaki, who
``became aware of impossibility to define'' for wild links what for PL links
is equivalent to the formal power series $\sum\beta^nz^n$ by a change of
variable.
\endabstract


\address Steklov Mathematical Institute, Division of Geometry and Topology;
\newline ul. Gubkina 8, Moscow 119991, Russia \endaddress
\curraddr University of Florida, Department of Mathematics;\newline
358 Little Hall, PO Box 118105, Gainesville, FL 32611-8105, U.S. \endcurraddr
\email melikhov\@math.ufl.edu, melikhov\@mi.ras.ru \endemail
\address Institute of Mathematics, Physics and Mechanics, University of
Ljubljana, 19 Jadranska cesta, 1001 Ljubljana, Slovenija \endaddress
\email dusan.repovs\@fmf.uni-lj.si \endemail
\endtopmatter

\document 

\head 1. Introduction \endhead

This paper can be read independently of {\it ``$n$-Quasi-isotopy I''}.
We relate $k$-quasi-isotopy to $k$-cobordism in \S3, which is entirely
independent of \S2 and only uses the following definition from \S1.
\S2 is concerned with $\omega$-quasi-isotopy.
This section illustrates various versions of $k$-quasi-isotopy by examples
and simple geometric observations, and reduces them to $k$-splitting and
$k$-contractibility.
Unless the contrary is explicitly stated, everything is understood to be in
the PL category.

\subhead \oldnos1.\oldnos1.\,Definition \endsubhead
Let $N$ be a compact $1$-manifold, $M$ an orientable $3$-manifold, and $k$
a nonnegative integer for purposes of this section and \S3, or an ordinal
number for purposes of \S2.
Two PL embeddings $L,L'\:(N,\partial N)\emb (M,\partial M)$ will be called
{\it (weakly) [strongly] $k$-quasi-isotopic} if they can be joined by
a generic PL homotopy where every singular level is a (weak) [strong]
$k$-quasi-embedding.
A PL map $f\:(N,\partial N)\to (M,\partial M)$ with precisely one double
point $x$ will be called a {\it (weak) [strong] $k$-quasi-embedding} if
there exist an arc $J_0$ with $\partial J_0=f^{-1}(x)$ and chains of
subpolyhedra $\{x\}=P_0\i\dots\i P_k\i M$ and $J_0\i\dots\i J_k\i N$
such that
\smallskip

(i) $f^{-1}(P_i)\i J_i$ for each $i\le k$;

(ii) $P_i\cup f(J_i)\i P_{i+1}$ for each $i<k$;

(iii) the inclusion $P_i\cup f(J_i)\emb P_{i+1}$ is null-homotopic (resp\.
induces zero homomorphisms on reduced integral homology) for each $i<k$;

(iv) $P_{i+1}$ is a compact PL $3$-manifold [resp\. closed PL $3$-ball]
and $J_{i+1}$ is a closed arc for each $i<k$.
\medskip

We could let the $P_i$'s with finite indices be arbitrary compact
subpolyhedra of $M$ (as we did in {\it ``$n$-Quasi-isotopy I''}); the above
situation is then restored by taking regular neighborhoods.
Also, by \cite{Sm3} or \cite{Ha}, condition (iii) can be weakened for
finite $i$ to

\smallskip
(iii$'$) the inclusion $P_i\cup f(J_i)\emb P_{i+1}$ is trivial
on $\pi_1$ (resp\. $H_1$) for finite $i<k$.
\bigskip

Let us consider small values of $k$ in more detail, in the case where
$M=S^3$ and $N$ is $mS^1:=S^1_1\sqcup\dots\sqcup S^1_m$.
Let $f\:mS^1\to S^3$ be a map with precisely one double point $f(p)=f(q)$.
Clearly, $f$ is a {\sl $0$-quasi-embedding}, in either of the $3$ versions,
iff it is a {\it link map}, i.e\. $f(S^1_i)\cap f(S^1_j)=\emptyset$ whenever
$i\ne j$.
Hence all $3$ versions of $0$-quasi-isotopy coincide with {\it link
homotopy}, i.e\. homotopy through link maps.
Note that if $f$ is a link map, its singular component splits into two
{\it lobes}, i.e\. the two loops with ends at $f(p)=f(q)$.

The map $f$ is a {\sl (weak) [strong] $1$-quasi-embedding} iff it is a link
map, and at least one of the two lobes of the singular component (namely,
$J_0$) is null-homotopic (resp\. null-homologous) [resp\. contained in a
PL $3$-ball] in the complement $X$ to the other components.
See \cite{MR; \S2} for some observations on $1$-quasi-isotopy.

Finally, $f$ is a {\sl $2$-quasi-embedding} iff it is a $1$-quasi-embedding,
so $f(J_0)$ is null-homotopic in $X$, and, for some arc $J_1\i S^1_i$
containing $f^{-1}(F(D^2))$, this null-homotopy $F\:D^2\to X$, which we
assume PL and generic, can be chosen so that every loop in
$F(D^2)\cup f(J_1)$ is null-homotopic in $X$.
Note that a point of the finite set $f^{-1}(F(\Int D^2))$ may be outside
$J_0$ as well as inside it (compare examples (i) and (iii) below).

\subhead \oldnos1.\oldnos2. Examples \endsubhead
The reader is encouraged to visualize

(i) a strong null-$(k-1)$-quasi-isotopy for the Milnor link $\M_k$ \cite{Mi2}
(for $k=4$ see Fig\. 1 ignoring the four disks), and for any of its twisted
versions (in Fig\. 1, imagine any number of half-twists along each disk);

(ii) a strong $k$-quasi-isotopy between any two twisted versions of $\M_k$
that only differ by some number of full twists in the rightmost clasp;

\def\epsfsize#1#2{0.75\hsize}
\fig 1
\def\epsfsize#1#2{\hsize}

(iii) a null-$(k-1)$-quasi-isotopy for the link $\W_k$, the $k$-fold
untwisted left handed Whitehead double of the Hopf link (for $k=3$ see
Fig\. 2a ignoring the three disks);

(iii$'$) a null-$(k-1)$-quasi-isotopy for any $k$-fold Whitehead double of
the Hopf link that is untwisted (i.e\. differs from $\W_k$ by at most one
positive half-twist) at all stages except possibly for the first or the last,
and is arbitrarily twisted at that stage (in Fig\. 2a, imagine any number of
half-twists along one of the two smaller disks and possibly one clockwise
180 degree rotation of the visible side of the other smaller disk and/or
the larger disk).

(iv) a weak null-$(k-1)$-quasi-isotopy for an arbitrarily twisted $k$-fold
Whitehead double of the Hopf link (any number of half-twists along each disk
in Fig\. 2a);

(v) a (weak) $k$-quasi-isotopy between any versions of $\W_k$, untwisted at
the first $k-1$ stages and arbitrarily twisted at the last stage (resp\.
arbitrarily twisted at all stages), that only differ from each other at the
last stage.

\fig 2

Note that it does not matter which component is doubled in an iterated left
handed untwisted Whitehead doubling of the Hopf link (and therefore $\W_k$ is
well-defined) by the symmetry of $\W_1$, i.e\. realizability of its
components' interchange by an ambient isotopy, applied inductively as shown
in Fig\. 2b.

\remark{Remark} We show in \S3 that $\M_k$ is not null-$k$-quasi-isotopic,
for each $k\in\N$.
\endremark

\proclaim{Conjecture 1.1} (a) $\W_k$ is not null-$k$-quasi-isotopic, even
weakly, for each $k\in\N$.

(b) $\W_k$ is not strongly null-$1$-quasi-isotopic, for each $k\in\N$.
\endproclaim

Certainly, the obvious null-$(k-1)$-quasi-isotopy in (iii) is not a strong
$1$-quasi-isotopy, since $\W_{k-1}$ (for $k\ge 1$) is not a split link
\cite{Wh} (see also \cite{BF; \S3}), and is not a weak $k$-quasi-isotopy,
since the components of $\W_{k-1}$ are $(k-1)$-linked \cite{Sm2}.

\subhead \oldnos1.\oldnos3. $n$-Linking, $n$-splitting and
$n$-contractibility\endsubhead
In 1937, S. Eilenberg term\-ed two disjoint knots $K_1\cup K_2\i S^3$
$0$-linked if they have nonzero linking number, and {\it $n$-linked} if
every subpolyhedron of $S^3\but K_2$, in which $K_1$ is null-homologous,
contains a knot, $(n-1)$-linked with $K_2$.
Note that if $K_1$ and $K_2$ are not $n$-linked, $K_1$ bounds a map of
a grope of depth $n+1$ into $S^3\but K_2$, i.e\. represents an element of
the $(n+1)^{\text{th}}$ derived subgroup of $\pi_1(S^3\but K_2)$.

Having proved that the components of $\Cal W_n$ are $n$-linked (as it had long
been expected), N. Smythe proposed in 1970 that the relation of being not
$n$-linked should be replaced by a stronger relation with substantially lower
quantifier complexity.
Compact subsets $A,B\i S^3$ are said to be {\it $n$-split}
\cite{Sm2; p\. 277}, if there exists a sequence of compact subpolyhedra
$A\i P_0\i\dots\i P_{n+1}\i S^3\but B$ such that each inclusion
$P_i\i P_{i+1}$ is trivial on reduced integral homology groups.
Thus ``$(-1)$-split'' means ``disjoint'', and two disjoint knots are
$0$-split iff they have zero linking number.
It is also not hard to see that the following are equivalent for a
$2$-component PL link $K_1\cup K_2$:
\roster
\item $K_1$ and $K_2$ are $1$-split;
\item $K_1$ and $K_2$ are not $1$-linked;
\item $K_1\cup K_2$ is a {\it boundary link}, i.e\. $K_1$ and $K_2$
bound disjoint Seifert surfaces.
\endroster
On the other hand, it turns out (see Theorem 2.8 below) that if $A$ and $B$
are compact subpolyhedra of $S^3$, there is an $n\in\N$ such that if $A$ and
$B$ are $n$-split, they are split by a PL embedded $S^2$.

An advantage of working in $S^3$ is that the relation of being $n$-split is
symmetric for compact subpolyhedra of $S^3$ \cite{Sm2}.
However, as in \cite{Sm3}, we can consider the same notion of $n$-splitting
for compact subsets $A,B\i M$ of an orientable $3$-manifold $M$.
Its interest for us stems from

\proclaim{Theorem 1.2} For each $n\in\N$, links
$L_1,L_2\:(I\sqcup\dots\sqcup I,\partial)\emb (M,\partial M)$ are (weakly)
$n$-quasi-isotopic if and only if they can be joined by a generic PL homotopy
$h_t$ such that for every $t$, each component of the image of $h_t$ is
$(n-1)$-contractible in the complement to the others (resp\. is $(n-1)$-split
from the union of others).
\endproclaim

The notion of $n$-contractibility was introduced recently by K. Kobayashi
\cite{Ko}.
We say that a compact subset $A\i M$ is {\it $n$-contractible} in $M$ if
there exists a sequence of compact subpolyhedra $A\i P_0\i\dots\i P_{n+1}\i M$
such that each inclusion $P_i\i P_{i+1}$ is null-homotopic.%
\footnote{In the definition of $n$-contractibility ($n$-splitting)
it is enough to require that each homomorphism $\pi_1(P_i)\to\pi_1(P_{i+1})$
(resp\. $H_1(P_i)\to H_1(P_{i+1})$) be zero, provided that $A$ is a
$1$-dimensional subpolyhedron of $M$ with no edges in $\partial M$
\cite{Ha}, \cite{Sm3}; compare \cite{Sm2; argument on p\. 278}.}

\demo{Proof} The `if' part follows immediately from the definitions by taking
regular neighborhoods.
The converse is also easy: just replace the manifolds $P_i$ from
the definition of (weak) $n$-quasi-embedding by their unions with the entire
singular component. \qed
\enddemo

This result does not hold for links in $S^3$ (see Remark (iii) below),
for which we could only find a faint version of it (see next proposition).
However, a link of $m$ circles in $S^3$ uniquely corresponds to a link
of $m$ arcs in the exterior of $m$ PL $3$-balls in $S^3$, and this bijection
descends to $n$-quasi-isotopy classes (as well as ambient isotopy classes
and PL isotopy classes).
This avoids the non-uniqueness of presentation of a link in $S^3$ as
the closure of a string link.

\subhead \oldnos1.\oldnos4. Brunnian links and Whitehead doubling \endsubhead
We now turn back to links in $S^3$.
Clearly, if $f$ is a weak $n$-quasi-embedding with a single double point, the
lobe $J_0$ of $f$ is $(n-1)$-split from the union of the non-singular
components.
Conversely, if the $i^{\text{th}}$ component of a link is null-homotopic in
a polyhedron, $(n-1)$-split from the union of the other components, then
the link is evidently weakly $n$-quasi-isotopic to a link with
the $i^{\text{th}}$ component split by a PL embedded $S^2$ from the other
components.
In particular:

\proclaim{Proposition 1.3}
If a component of a Brunnian link is $(n-1)$-split from the union of the
other components, then the Whitehead doubling of the link on this component
(with arbitrary twisting) is weakly null-$n$-quasi-isotopic.
\endproclaim

\remark{Remarks. (i)}
As shown by the Borromean rings, Proposition 1.3 cannot be ``desuspended'',
that is, a Brunnian link with a component $n$-split from the union of
the other components is not necessarily weakly null-$n$-quasi-isotopic.
\endremark

\remark{(ii)} Assuming Conjecture 1.1, $n-1$ cannot be replaced with $n-2$
in Proposition 1.3, since the components of $\W_{n-1}$ are $(n-2)$-split.
\endremark

\remark{(iii)}
Neither version of null-$n$-quasi-isotopy implies even $1$-splitting, since
$\M_{k+1}$ is strongly null-$n$-quasi-isotopic, but is not a boundary link,
as detected by Cochran's invariants (cf\. \S3).
\endremark
\medskip

Once Whitehead doubling appeared above, it is tempting to study the behavior
of our filtration under this operation.
It is easy to see that if two links are weakly $n$-quasi-isotopic with support
in the $i^{\text{th}}$ component, then the links, obtained by Whitehead
doubling (with arbitrary twisting) on all components except for the
$i^{\text{th}}$ one, are weakly $(n+1)$-quasi-isotopic with support in
the $i^{\text{th}}$ component.
Two possible ways to deal with the Whitehead double $K'$ of the $i^{\text{th}}$
component $K$ (drag $K'$ along the given homotopy of $K$ or unclasp $K'$ in
a regular neighborhood of $K$) result in the following two corollaries.

\proclaim{Proposition 1.4}
(a) Weakly $1$-quasi-isotopic links have weakly $2$-quasi-isotopic
Whitehead doubles.
Weakly $2$-quasi-isotopic links with trivial pairwise linking numbers have
weakly $3$-quasi-isotopic Whitehead doubles.

(b) The Whitehead double of the Borromean rings is weakly
null-$2$-quasi-isotopic.
\endproclaim

Here the Whitehead doubling is performed on all components, with an arbitrary
(but fixed) twisting on each.
Part (a) does not generalize, since the property of being a boundary
link is not preserved under $3$-quasi-isotopy (see Remark (iii) above).

\demo{Proof} Let us verify e.g\. the second part of (a).
If $f\:mS^1\to S^3$ is a singular link in a generic PL weak $2$-quasi-isotopy
between the given links, let $P_1$ be as in the definition, and $N\sqcup N'$
be a regular neighborhood of $f(mS^1)$ in $S^3$ such that $f^{-1}(N)$ is
the singular component, say $S^1_i$, and $N'\cap P=\emptyset$.
We can assume that $P_1\cap N$ is connected, and that $N$ (resp\. $N'$)
contains the Whitehead doubles of the two resolutions $K_+,K_-$ of $f(S^1_i)$
(resp\. $Wh(f(S^1_j))$ for all $j\ne i$).
Then each of the four singular knots in the obvious homotopy between $Wh(K_+)$
and $Wh(K_-)$ is a $1$-quasi-embedding, if regarded as a map into $N\cup P_1$.
The hypotheses imply that every cycle in either $P_1$ or $N$ is null-homologous
in $S^3\but N'$, hence so is every cycle in $P_1\cup N$.
Finally, every cycle in $S^3\but N'$ is null-homologous in
$S^3\but\bigcup_{j\ne i} Wh(S^1_j)$. \qed
\enddemo

\head 2. $\omega$-Quasi-isotopy \endhead

In this section we study how the relations of $k$-quasi-isotopy approximate
the relation of {\it PL isotopy}, i.e\. PL homotopy through embeddings.

\subhead \oldnos2.\oldnos1. Basic observations \endsubhead
First of all we get rid of infinite ordinals in the definition of
$k$-quasi-isotopy for $k=\omega,\omega+1$ and $\omega_1$.

\proclaim{Proposition 2.1} Let $f\:N\to M$, its double point $\{x\}=P_0$,
the $3$-manifolds $P_1\i P_2\i\dots$ and the arcs $J_0\i J_1\i\dots$ be as in
the definition of (weak) [strong] $k$-quasi-embedding in \S1.
Then

(a) $f$ is a (weak) [strong] $\omega$-quasi-embedding iff the $P_i$, $J_i$'s
exist for all finite $i$.

(b) $f$ is a (weak) $(\omega+1)$-quasi-embedding iff the $P_i$, $J_i$'s exist
for all finite $i$ and $\bigcup_{i<\omega}P_i$ is contained in a compact
$3$-manifold $Q\i M$.

(c) $f$ is a (weak) [strong] $\omega_1$-quasi-embedding iff $x$ is contained
in a contractible (acyclic) compact $3$-manifold [closed $3$-ball] $P_*$
such that $f^{-1}(P_*)$ is an arc $J_*$.

(d) Two links $L,L'\:(N,\partial N)\emb (M,\partial M)$ are strongly
$\omega_1$-quasi-isotopic iff they are PL isotopic.
\endproclaim

By (d), PL isotopy implies either version of $k$-quasi-isotopy for any $k$;
conversely, each $k$-quasi-isotopy is supposed to give an approximation of
PL isotopy.
This explains why the $P_i$'s were required to be compact in \S1 when $i=i'+1$
for some $i'$; the point of our dropping this requirement when this is not
the case (i.e\. when $i$ is a {\it limit} ordinal) is that the interesting
relation of $\omega$-quasi-isotopy would otherwise fall out of the hierarchy
due to the resulting shift of infinite indices.%
\footnote{The extra care about limit ordinals $>\omega$ will turn out to be
superfluous in light of the following results.
Note, however, that the $k$-quasi-concordance filtration on classical links
\cite{MR; Remark (ii) at the end of \S1} may well be highly nontrivial when
similarly extended for infinite ordinals.}

\demo{Proof. (a)} Set $P_\omega=\bigcup_{i<\omega}P_i$ and
$J_\omega=\bigcup_{i<\omega}J_i$.
Then the only condition $f^{-1}(P_\omega)\i J_\omega$ imposed on $P_\omega$
and $J_\omega$ holds due to $f^{-1}(P_i)\i J_i$ for finite $i$.
\enddemo

\demo{(b)} Set tentatively $P_\omega=\bigcup_{i<\omega} P_i$ and
$J_\omega=\bigcup_{i<\omega} J_i$.
If the closure of $J_\omega$ happens to be the entire singular component
$S^1_i$, we puncture $S^1_i$ at the singleton $pt=S^1_i\but J_\omega$ and
$M$ at $f(pt)$, and sew in a trivial closed PL ball pair $(B^3,B^1)$.
(This in effect redefines $J_k,\dots,J_\omega$ and $P_k,\dots,P_\omega$
starting from some $k$, similarly to the explicit redefinition that follows.)
Next, consider a map $F\:M\to M$ shrinking a regular neighborhood $R$ of
the $1$-manifold $Z:=\Cl{f(N\but J_1)}$ rel $\partial J_1$ onto $Z$ so that
all non-degenerate point inverses are closed PL $2$-disks, and sending
$M\but R$ homeomorphically onto $M\but Z$.
Redefine each $P_i$ to be $F^{-1}(P_i)$ for $2\le i\le\omega$; it is easy
to see that the conditions (i)--(iv) from \S1 still hold for these.
Set $P_{\omega+1}=\Cl{F^{-1}(Q)\but F^{-1}(Z')}$, where
$Z'=f(N\but J_\omega)$.
This is a compact $3$-manifold (since $Q$ is), and
$J_{\omega+1}:=f^{-1}(P_{\omega+1})=\Cl{J_\omega}$ is a closed arc.
Moreover, the inclusion $P_\omega\cup f(J_\omega)\i P_{\omega+1}$ holds
(since $f^{-1}(P_\omega)=J_\omega$ due to conditions (i) and (ii) for finite
$i$, and since $P_\omega\i Q$), and is null-homotopic/null-homologous (since
$P_\omega$ is contractible/acyclic due to condition (iii) for finite $i$).
\qed
\enddemo

\demo{(c)} The $3$-manifolds $P_i$ and the arcs $J_i$ have to stabilize
at some countable stage since $M$ and $N$ are separable. \qed
\enddemo

\demo{(d)} Let $h_s$ be a singular level in a strong $\omega$-quasi-isotopy
$h_t\:(N,\partial N)\to (M,\partial M)$, and let the $3$-ball $P_*$ and
the arc $J_*$ be as in (c) corresponding to $f=h_s$.
Without loss of generality, there is an $\eps>0$ such that for each
$t\in [s-\eps,s+\eps]$, the level $h_t$ is an embedding, coinciding with
$h_s$ outside $J_*$.
The transition between the links $h_{s-\eps}$ and $h_{s+\eps}$ can now be
realized by a PL isotopy that is conewise on $J_*$ and first shrinks
the local knot of $h_{s-\eps}$ occurring in $P_*$ to a point and then inserts
the local knot of $h_{s+\eps}$ occurring in $P_*$ by an inverse process. \qed
\enddemo

\proclaim{Lemma 2.2} \cite{Sm1} A knot in a $3$-manifold $M$, PL isotopic
to the unknot, is contained in a ball in $M$.
\endproclaim

Under the unknot in a $3$-manifold $M$ we understand the ambient isotopy
class in $M$ of the unknot in some ball in $M$ (since all balls in $M$ are
ambient isotopic).
A PL isotopy may first create a local knot $K$, then push, say, a homology
ball ``through a hole'' in $K$, and finally shrink $K$ back to a point,
so Lemma 2.2 is not obvious.
In fact, its higher-dimensional analogue is false (for PL isotopy of $S^n$ in
$S^{n+2}\but S^n$) \cite{Ro}.
There also exist two PL links $L,L'\:2S^1\emb S^1\x S^2$ that are PL isotopic
but not ambient isotopic, even though each component of $L$ is ambient
isotopic with the corresponding component of $L'$ \cite{Ro}.

For convenience of the reader we shall give an alternative proof of
Lemma 2.2, based on an approach different from Smythe's.

\demo{Proof} This is in the spirit of the uniqueness of knot factorization.
Suppose that a knot $K'\i\Int M$ is obtained from a knot $K\i\Int M$
by a PL isotopy with support in a ball $B\i\Int M$, which intersects $K$
(hence also $K'$) in an arc.
Assuming that $K$ is contained in a ball, we show that so is $B\cup K$,
hence also $K'$.
We may view ball as a special case of punctured ball, i.e\. a homeomorph
of the exterior in $S^3$ of a nonempty collection of disjoint balls.
Let $n$ be the minimal number of components in the intersection of
$\partial B$ with the boundary of a punctured ball $P\i\Int M$,
containing $K$.

We claim that $n=0$.
Let $C$ be a circle in $\partial B\cap\partial P$ that is innermost on
$\partial B$ (with $\partial B\cap K$ considered ``outside'').
Then $C$ bounds a disk $D\i\partial B$ with $D\cap\partial P=C$.
If $D$ lies outside $P$, we attach to $P$ the embedded $2$-handle with core
$D$, thus increasing the number of components in $\partial P$, but decreasing
the number of components in $\partial B\cap\partial P$.
If $D$ is inside $P$, we subtract the $2$-handle from $P$, which by
Alexander's Sch\"onflies Theorem splits $P$ into two punctured balls, and
discard the one which is disjoint from $K$.

Thus we may assume that each component of $\partial P$ is disjoint form
$\partial B$.
By Alexander's Theorem the components of $\partial P$, contained in $B$,
bound balls in $B$.
Since $B$ contains only an arc of $K$, $P$ itself cannot be contained in $B$,
so these balls are disjoint from $\Int P$.
Adding them to $P$ makes $P$ have no boundary components in $B$, hence contain
$B$.
If this makes $P$ have no boundary whatsoever, $M$ must be $S^3$, and there
is nothing to prove then.
Otherwise we may connect the components of $\partial P$ by drilling holes
in $P\but (B\cup K)$, so as to obtain a ball containing $B\cup K$. \qed
\enddemo

\remark{Remark}
D. Rolfsen asked \cite{Ro} whether every two-component link in a $3$-manifold
that is isotopic to the unlink is a split link (i.e.\ a link whose components
are contained in disjoint balls).
In fact, this had been established earlier by Smythe \cite{Sm1}.
The above argument yields a different proof of this fact.
Indeed, if $L$ is isotopic to a split link $L'$ by a PL isotopy with support
in the $i^{\text{\rm th}}$ component, by the proof of Lemma 2.2, the
$i^{\text{\rm th}}$ component of $L$ is contained in a ball $B$ disjoint from
the other components.
Although $B$ may intersect the disjoint balls containing the other components
by the hypothesis, this can be remedied via the same argument.
\endremark

\subhead \oldnos2.\oldnos2. Simply-connected $3$-manifolds and
finiteness results \endsubhead
Recall the construction of the Whitehead contractible open manifold $W$.
Consider a nested sequence of solid tori $\dots T_1\i T_0\i S^3$, where each
$T_{i+1}\cup S^3\but T_{i+1}$ is equivalent, by a homeomorphism of $S^3$,
to a regular neighborhood of the Whitehead link $\W_1$ in $S^3$.
Then $W$ is the union of the ascending chain of solid tori
$S^3\but T_0\i S^3\but T_1\i\dots$.
Now $W$ is not homeomorphic to $\R^3$, since a meridian
$K_0:=\partial D^2\x\{pt\}$ of $T_0\cong D^2\x S^1$ is not contained in any
(piecewise-linear) $3$-ball in $W$; otherwise the ball would be disjoint
from some $T_n$, hence $\W_n$ (see Fig\. 2 above) would be a split link,
which is not the case \cite{Wh} (see also \cite{BF; \S3}).

The preceding lemma implies that $K_0$ is not PL isotopic to the unknot
in $W$.
On the other hand, since the Whitehead manifold $W$ is contractible, it is
not hard to see that any knot $K$ in $W$ is $\omega$-quasi-isotopic to
the unknot.
Indeed, if $f$ is a singular level in a generic PL homotopy between $K$ and
the unknot, we set $J_0$ to be any of the two lobes, and construct
the polyhedra $P_j$ and the arcs $J_j$, $1\le j<\omega$, as follows.
Assuming that $P_j$ is a handlebody, we set $P_{j+1}$ to be a regular
neighborhood of $P_j$ union the track of a generic PL null-homotopy of
a wedge of circles onto which $P_j\cup f(J_j)$ collapses.
By \cite{Sm3} or \cite{Ha}, the null-homotopy can be chosen so that $P_{j+1}$
is again a handlebody.
Since the track is $2$-dimensional and generic, $f^{-1}(P_{j+1})$ is not
the whole circle, and we set $J_{j+1}$ to be the smallest arc containing
$f^{-1}(P_{j+1})$.
Then $f$ is an $\omega$-quasi-embedding by Proposition 2.1(a).

In fact, any contractible open $3$-manifold other than $\R^3$ contains
a knot which is not contained in any $3$-ball \cite{CDG}.
Thus the above argument proves

\proclaim{Proposition 2.3} Any contractible open $3$-manifold other than
$\R^3$ contains a knot, $\omega$-quasi-isotopic but not PL isotopic to
the unknot.
\endproclaim

By Bing's characterization of $S^3$ (see \cite{Ei} for a proof and references
to 4 other proofs), any closed $3$-manifold $M$ other than $S^3$ contains
a knot $K$ that is not contained in any ball, hence is not PL isotopic to
the unknot.
If $M$ is a homotopy sphere, clearly $K$ is $\omega_1$-quasi-isotopic
to the unknot.
(Indeed, let $f$ be a singular level in a generic PL homotopy between $K$
and the unknot, and let $B$ be a PL $3$-ball in $M$ meeting $f(mS^1)$ in
an arc; then $P_*:=\Cl{M\but B}$ is contractible and $f^{-1}(P_*)$ is an
arc.)
This proves

\proclaim{Proposition 2.4} Any homotopy $3$-sphere other than $S^3$ (if
exists) contains a knot, $\omega_1$-quasi-isotopic but not PL isotopic to
the unknot.
\endproclaim

Let us now turn to positive results.

\proclaim{Theorem 2.5} Two links in a compact orientable $3$-manifold are
weakly $\omega$-quasi-isotopic if and only if they are weakly
$\omega_1$-quasi-isotopic.
\endproclaim

This is in contrast with the existence of compact orientable $3$-manifolds
such that the lower central series of the fundamental group does not
stabilize until the $2\omega^{\text{th}}$ stage \cite{CO} (see also
\cite{Mih}).
However, it is an open problem, to the best of our knowledge, whether
a compact $3$-manifold group or indeed a finitely presented group may have
derived series of length $>\omega$.
(Theorem 2.8 below implies that embedded loops in a compact $3$-manifold $M$
that bound {\it embedded} gropes of depth exceeding certain finite number
$n(M)$ are contained in the intersection of the transfinite derived series of
$\pi_1(M)$.)

\proclaim{Corollary 2.6} Let $M$ be a compact orientable $3$-manifold
containing no nontrivial homology balls.
Then the relations of $\omega$-quasi-isotopy, weak $\omega$-quasi-isotopy,
strong $\omega$-quasi-isotopy and PL isotopy coincide for links in $M$.
\endproclaim

Note that by Alexander's Sch\"onflies Theorem any $3$-manifold, embeddable in
$S^3$, contains no nontrivial homology balls.

\proclaim{Corollary 2.7} The assertions of Theorem 2.5 and Corollary 2.6 hold
for links in a non-compact orientable $3$-manifold if $\omega$ is replaced by
$\omega+1$.
\endproclaim

The first of the two corollaries above follows by part (d) and the second
by part (a) along with the (trivial) implication ``only if'' of part (b) of
Proposition 2.1.

\demo{Proof of Theorem 2.5} Let $Q^3$ be the manifold and $f\:mS^1\to Q$
a weak $\omega$-quasi-embedding.
Let $S^1_i$ be the singular component, and $J_0$ be as in the definition of
$k$-quasi-embedding.
Let $D^3$ be a regular neighborhood of some $y\in f(S^1_i\but J_0)$ in
$N^3$, meeting the image of $f$ in an arc $D^1$.
Then $A:=\Cl{f(S^1_i)\but D^1}$ is $n$-split from $B:=f(mS^1\but S^1_i)$
in $M:=\Cl{Q^3\but D^3}$ for every $n\in\N$ (cf\. Theorem 1.2).
Thus the assertion follows from
\enddemo

\proclaim{Theorem 2.8} Let $M$ be a compact connected orientable
$3$-manifold and $A,B$ nonempty compact subpolyhedra of $M$.
Then there exists an $n\in\N$ such that if $A$ is $n$-split from $B$, then
$A$ is contained in a PL homology ball in $M\but B$.
\endproclaim

In the case $M=S^3$ this was claimed in \cite{Sm2}, but the proof there
appears to be valid only if both $A$ and $B$ are connected (Lemma on
p\. 279 is incorrect if $A$ is disconnected, due to a misquotation on
p\. 280 of Theorem from p\. 278, and no proof is given when $B$ is
disconnected).
The idea of our proof is close to Smythe's, but essential technical changes
are needed in case $A$ or $B$ is disconnected, $\partial M\ne\emptyset$ or
$H_1(M)\ne 0$.
Interplay of these cases brings additional complications: it is certainly
unnecessary to invoke Ramsey's Theorem if either $B$ is connected or
$\partial M=\emptyset$ (argue instead that for each $s$ there is an $i$
such that $H_1(M^s_i)=0$).
Note that in proving Theorem 2.5 for $3$-component links in $S^3$ one
already encounters both nonempty $\partial M$ and disconnected $B$ in
Theorem 2.8.

Note that the knot $K_0$ in the Whitehead manifold $W$ (see above) and
the Whitehead continuum $S^3\but W$ are mutually $n$-split in $S^3$ for all
$n\in\N$, but no PL embedded $S^2$ splits them, since each $\W_n$ is not
a split link.
Another such example is based on Milnor's wild link $\Cal M_\infty$
\cite{Mi2; p\. 303}, which can be approximated by the links $\Cal M_n$
(see Fig\. 1 above).
Remove from $S^3$ a PL $3$-ball, disjoint from the tame component and meeting
the wild component in a tame arc, and let
$\mu\:(S^1\sqcup I,\partial I)\emb (B^3,\partial B^3)$ denote the resulting
wild link.
Then $\mu(I)$ is $n$-split from $\mu(S^1)$ in $B^3$ for all $n\in\N$, but is
not contained in any PL $3$-ball in $B^3\but\mu(S^1)$, since each $\Cal M_n$ is
not PL isotopic to the unlink (cf\. \S3).

\subhead \oldnos2.\oldnos3. Proof of the finiteness results \endsubhead
The proof of Theorem 2.8 is based on the Kneser--Haken Finiteness
Theorem, asserting that there cannot be infinitely many disjoint
incompressible, $\partial$-incompres\-sible pairwise non-parallel surfaces%
\footnote{A connected $2$-manifold $F$, properly PL-embedded in an orientable
$3$-manifold $M$, is called {\it compressible} if either $F=\partial B$ for
some $3$-ball $B\i M$ or there exists a $2$-disk $D\i M$ such that
$D\cap F=\partial D$ and $D\cap F\ne\partial D'$ for any $2$-disk $D'\i F$.
Next, $F$ is {\it $\partial$-compressible} if either
$F=\Cl{\partial B\but\partial M}$ for some $3$-ball $B\i M$, meeting
$\partial M$ in a disk, or there exists a $2$-disk $D\i M$, meeting
$\partial M$ in an arc and such that
$D\cap F=\Cl{\partial D\but\partial M}$ and
$D\cap F\ne\Cl{\partial D'\but\partial M}$ for any $2$-disk $D'\i M$.
Finally, $F$ and $F'$ are {\it parallel} in $M$ if $(F',\partial F')$ is
the other end of a collar $h(F\x I,\partial F\x I)$ of
$(F,\partial F)=h(F\x\{0\},\partial F\x\{0\})$ in $(M,\partial M)$.}
in a compact $3$-manifold $M$ \cite{Go}, \cite{Ja; III.20}.
(Erroneous theorem III.24 in \cite{Ja}, which was concerned with weakening
the condition of $\partial$-incompressibility, is corrected in \cite{FF}).
More precisely, we will need a slightly stronger version of this result,
where the manifold $M$ is obtained by removing a compact surface with
boundary from the boundary of a compact $3$-manifold $\hat M$.
(A proper surface in $M$, $\partial$-incompressible in $M$, need not be
$\partial$-incompressible in $\hat M$.)
This strengthening is proved in \cite{Ma; 6.3.10} in the case where $M$
is irreducible and $\partial$-irreducible; the general case can be reduced
to this case by the same argument \cite{Ja}.

\definition{Notation}
Fix an $n\in\N$, and let $M_0$ and $\Cl{M\but M_{n+1}}$ be disjoint regular
neighborhoods of $A$ and $B$ in $M$.
A collection $(M_1,\dots,M_n)$ of $3$-dimensional compact PL submanifolds of
$M$ such that $M_i\cap\Cl{M\but M_{i+1}}=\emptyset$ and
$i_*\:H_1(M_i)\to H_1(M_{i+1})$ is zero for each $i=0,\dots,n$ will
be called a {\it pseudo-splitting}.
Its {\it complexity} is defined to be
$$c(M_1,\dots,M_n):=\sum_{i=1}^n\sum_{\quad F\in\pi_0(\Fr M_i)}
\left(\vphantom{\hat h}\rk H_1(F)+\rk H_{-1}(\partial F)\right)^2,$$
where $H_{-1}(X)=0$ or $\Z$ according as $X$ is empty or not, and $\pi_0$
stands for the set of connected components.
\enddefinition

\proclaim{Lemma 2.9} {\rm (compare \cite{Mc1; p\. 130})}
If $(M_1,\dots,M_n)$ is a pseudo-splitting of minimal complexity, each
component of $\Fr M_1\cup\dots\cup\Fr M_n$, which is not a sphere or a disk,
is incompressible and $\partial$-incompressible in $M\but (A\cup B)$.
\endproclaim

\demo{Proof} Suppose that a component of some $\Fr M_i$, $1\le i\le n$, which
is not a sphere or a disk, is ($\partial$-)compressible in $M\but (A\cup B)$.
By the innermost circle argument, a component $F$ of some $\Fr M_j$,
$1\le j\le n$, which is not a sphere or a disk, is compressible (resp\.
compressible or $\partial$-compressible) in
$M\but(M_{j-1}\cup\Cl{M\but M_{j+1}})$.
If the ($\partial$-)compressing disk $D$ lies outside $M_j$, attach to $M_j$
an embedded $2$-handle with core $D$ (resp\. cancelling $1$-handle and
$2$-handle with cores $D\cap\partial M$ and $D$) and denote the result
by $M'_j$.
By Mayer--Vietoris, $i_*\:H_1(M_j)\to H_1(M'_j)$ is epic.
If $D$ lies inside $M_j$, remove the handle from $M_j$, and denote the result
again by $M'_j$; by Mayer--Vietoris, $i_*\:H_1(M'_j)\to H_1(M_j)$ is monic.
In any case, it follows that $(M_1,\dots,M'_j,\dots,M_n)$ is
a pseudo-splitting.

If $\partial D$ (resp\. $\Cl{\partial D\but\partial M}$) is not
null-homologous in $(F,\partial F)$, the result $F'\i\Fr M'_j$ of the surgery
on $F$ is connected, and $\rk H_1(F')<\rk H_1(F)$.
Since $\partial F'=\emptyset$ iff $\partial F=\emptyset$, we have that
$c(M_1,\dots,M_j',\dots,M_n)<c(M_1,\dots,M_n)$ in this case.
In the other case $F'$ consists of two connected components $F'_+$ and $F'_-$
such that $\rk H_1(F'_+)+\rk H_1(F'_-)\le\rk H_1(F)$.
If one of the summands, say $\rk H_1(F'_+)$, is zero, we must be in the case
of non-$\partial$ compression, and $F'_+$ must be a disk.
Hence if $c_F$ denotes $\rk H_1(F)+\rk H_{-1}(\partial F)$, both $c_{F'_+}$
and $c_{F'_-}$ are always nonzero, whereas $c_{F'_+}+c_{F'_+}\le c_F+1$.
Thus again $c(M_1,\dots,M_j',\dots,M_n)<c(M_1,\dots,M_n)$. \qed
\enddemo

\proclaim{Lemma 2.10} {\rm (compare \cite{Sm2; p\. 278})}
Let $(M_1,\dots,M_n)$ be a pseudo-splitting.

(a) If a component $N$ of $M_i$ is disjoint from $M_{i-1}$,
$(M_1,\dots,M_i\but N,\dots,M_n)$ is a pseudo-splitting.
If a component $N$ of $\Cl{M\but M_i}$ is disjoint from $\Cl{M\but M_{i+1}}$
and either $i<n$ or $H_1(M)\to H_1(N,\Fr N)$ is zero,
$(M_1,\dots,M_i\cup N,\dots,M_n)$ is a pseudo-splitting.

(b) There exists a pseudo-splitting $(M_1^*,\dots,M_n^*)$ such that

\noindent
{\rm (i)} each $\Fr M_i^*\i\Fr M_i$;

\noindent
{\rm (ii)} each $i_*\:H_0(M_{i-1}^*)\to H_0(M_i^*)$ and
$i_*\:H_0(M\but M_{i+1}^*)\to H_0(M\but M_i^*)$ are onto;

\noindent
{\rm (iii)} no component of $M_1^*$ is disjoint from $A$, and
$H_1(M)\to H_1(N,\partial N)$ is nonzero for each component $N$ of
$\Cl{M\but M_n^*}$, disjoint from $B$.
\endproclaim

\demo{Proof}
First note that (b) follows by an inductive application of (a).
The first assertion of (a) is obvious; we prove the second.
Assume that $H_1(M)\to H_1(N,\Fr N)$ is zero.
$H_1(M_i\cup N)\to H_1(N,\Fr N)$ factors through $H_1(M)$, hence is zero.
Then $i_*\:H_1(M_i)\to H_1(M_i\cup N)$ is onto, and the assertion follows.

If $i<n$, $H_1(M_{i+1})\to H_1(N,\Fr N)$ factors through $H_1(M_{i+2})$,
hence is zero.
On the other hand,
$H_1(M_{i+1})\to H_1(M_{i+1},M_i)\simeq H_1(N,\partial N)\oplus
H_1(M_{i+1},M_i\cup N)$ is monic.
So $H_1(M_{i+1})\to H_1(M_{i+1},M_i\cup N)$ is monic and
$H_1(M_i\cup N)\to H_1(M_{i+1})$ is zero. \qed
\enddemo

\demo{Proof of Theorem 2.8}
Let $D$ be $\rk H_1(M)$ plus the number of elementary divisors of
$\Tors H_1(M)$, and set $S=\rk H_0(A)+\rk H_0(B)+D+\rk H_1(M)$.
Recall the simplest case of Ramsey's Theorem: for each $k,l\in\N$ there
exists an $R(k,l)\in\N$ such that among any $R(k,l)$ surfaces in a
$3$-manifold either some $k$ are pairwise parallel, or some $l$ are pairwise
non-parallel.
Let $h$ be the number given by the Haken Finiteness Theorem for
the $3$-manifold $M\but (A\cup B)$.
Set $r_0=2$, $r_{i+1}=R(r_i,h)$ (so $r_1=h+1$) and finally $n=Sr_S$.

Since $A$ is $n$-split from $B$, there exists a pseudo-splitting
$(M'_1,\dots,M'_n)$ of minimal complexity.
Feed it into Lemma 2.10(b), and let $(M_1,\dots,M_n)$ denote the resulting
pseudo-splitting.
By Lemma 2.9 and (i) of Lemma 2.10(b), each component of
$\Fr M_1\cup\dots\cup\Fr M_n$, which is not a sphere or a disk, is
incompressible and $\partial$-incompressible in $M\but (A\cup B)$.
By (ii) and (iii) of Lemma 2.10(b), the same holds for sphere and disk
components.

If $N_1,\dots,N_d$ are the components of $\Cl{M\but M_n}$, disjoint from $B$,
the image $I_i$ of $H_1(M)$ in each $H_1(N_i,\Fr N_i)$ is nontrivial.
$H_1(M)$ surjects onto the subgroup $\oplus I_i$ of
$H_1(\bigcup N_i,\bigcup\partial N_i)$, so $d\le D$.
On the other hand, $i_*\:H_1(M_i)\to H_1(M)$ is zero for each $i\le n$,
and the image of $H_1(M\but M_i)$ in $H_1(M)$ contains that of
$H_1(M\but M_{i+1})$.
Hence the Mayer--Vietoris image of $H_1(M)$ in $H_0(\Fr M_i)$ changes at most
$\rk H_1(M)$ times as $i$ runs from $1$ to $n$.
Since $n=Sr_S$, by the pigeonhole principle we can find an $i_0$ such that
whenever $i_0\le i<i_0+r_S$, each component of $M_{i+1}\but M_i$ is adjacent
to the unique component of $\Fr M_i$ and also to the unique component of
$\Fr M_{i+1}$.

Let $\Fr M_i^s$ denote the component of $\Fr M_i$, corresponding to
$s\in\{1,\dots,S\}$ under the epimorphism
$H_0(A)\oplus H_0(B\cup\bigcup N_i)\oplus H_1(M)\to H_0(\Fr M_i)$.
Since $r_S=R(r_{S-1},h)$, by Ramsey's Theorem either some $r_{S-1}$ of
the surfaces $\Fr M^S_i$, $i_0\le i\le i_0+r_S$, are pairwise parallel in
$M\but (A\cup B)$, or some $h$ are pairwise non-parallel; but the latter
cannot be by Haken's Theorem.
Proceeding by induction, we eventually boil down to $r_0=2$ indices $i,j$
such that $i_0\le i<j\le i_0+r_S$ and $\Fr M^s_i$, $\Fr M^s_j$ are parallel
for each $s$.
Since each component of $M_j\but M_i$ is $M_j^s\but M_i^s$ for some $s$,
we obtain that $\Fr M_i$ and $\Fr M_j$ are parallel in $M\but (A\cup B)$.
Then $i_*\:H_1(M_i)\to H_1(M_j)$ is an isomorphism.
Since $(M_1,\dots,M_n)$ is a pseudo-splitting, it also is the zero map, so
$H_1(M_i)=0$.
Then $H_1(\partial M_i)$ is a quotient of $H_2(M_i,\partial M_i)=H^1(M_i)=0$,
hence $\partial M_i$ is a disjoint union of spheres.

Since $A$ is $0$-split from $B$, all components of $A$ are contained in
the same component of $M\but B$, and by Lemma 2.11 below all components of
$B\cup\partial M\but A$ are contained in the same component of $M\but A$.
If a component of $\partial M$ were contained in $A$, triviality of
$i_*\:H_2(A)\to H_2(M\but B)$ would imply $B=\emptyset$, contradicting the
hypothesis.
Hence the components of $\partial M_i$ can be connected by tubes in
$M\but (A\cup B)$, first by subtracting $1$-handles avoiding $A$ from
the components of $M_i$ and then by joining the amended components of $M_i$
by $1$-handles avoiding $B$.
The resulting $3$-manifold $M_i^+$ is connected, $\partial M_i^+$ is
a sphere, and clearly $H_1(M_i^+)=H_1(M_i)=0$, hence $M_i^+$ is a homology
ball. \qed
\enddemo

\proclaim{Lemma 2.11} If $A,B$ are disjoint subpolyhedra of a compact
connected orientable $3$-manifold $M$ and
$H_2(A)@>\sssize i_*>>H_2(M\but B)$ is zero,
$\tilde H^0(M\but A)@>\sssize i_*>>\tilde H^0(B\cup\partial M\but A)$
is zero.
\endproclaim

\demo{Proof}
By the hypothesis, $H_2(M\but B)\to H_2(M\but B,A\cap\partial M)$ is monic.
Hence $i_*\:H_2(A,A\cap\partial M)\to H_2(M\but B,A\cap\partial M)$ is zero,
and by Lefschetz duality so is
$H^1(M,M\but A)\to H^1(M,B\cup\partial M\but A)$.
On the other hand, $\tilde H^0(M)=0$, so
$\delta^*\:\tilde H^0(B\cup\partial M\but A)\to H^1(M,B\cup\partial M\but A)$
is monic, and the assertion follows. \qed
\enddemo

\subhead \oldnos2.\oldnos4. Acyclic open sets and strong
$\omega$-quasi-isotopy \endsubhead
There is another way to view Theorem 2.5 and its proof.
Showing that $\omega$-quasi-isotopy implies PL isotopy for links in a compact
$3$-manifold $Q$ amounts%
\footnote{Indeed, let $f\:mS^1\to Q$ be an $\omega$-quasi-embedding.
By taking increasingly loose regular neighborhoods in $Q$, we may assume
without loss of generality that $P_i$ lies in the interior of $P_{i+1}$
in $Q$ for each $i<\omega$.
Then $\bigcup_{i<\omega} P_i$ is contractible and open in $Q$, and this
property will preserve under the construction from the proof of
Proposition 2.1(b).
Using the final notation from that proof, $U:=P_\omega\cup(D\but\partial D)$,
where $D$ denotes the pair of disks $F^{-1}(\partial J_{\omega+1})$, is
a contractible open subset of $M:=P_{\omega+1}$ and contains
$J:=f(J_{\omega+1})$.}
to replacing a contractible open neighborhood $U$ of an arc $J$, properly
PL immersed with one double point into a compact $3$-manifold $M$, by
a closed PL ball neighborhood $V$.
The natural absolute version of this problem is to do the same for a PL knot
$K\i\Int M$.
Note that the latter is clearly impossible if one wishes to additionally
demand either $V\i U$ (by letting $U$ be the Whitehead manifold $W$ and $K$
be the meridian $K_0$ of $T_0$) or $V\supset U$
(by letting $U$ be Alexander's horned ball in $S^3$ and $M$ be
$\Cl{S^3\but T}$, where $T$ is an essential solid torus in $S^3\but\Cl{U}$).
However, it is possible to do this without such additional restrictions if
$M$ embeds in $S^3$, since Theorem 2.8 (with $B=$ point) implies

\proclaim{Corollary 2.12} If $U$ is an acyclic open subset of a compact
orientable $3$-manifold $M$, every compact subset of $U$ is contained in
a PL-embedded homology $3$-ball.
\endproclaim

This was proved by McMillan in the case where $M$ is irreducible,
$\partial M\ne\emptyset$ and $U\i\Int M$ \cite{Mc2; Theorem 2 + Remark 1}
(see also \cite{MT; Proposition} for the above formulation and
\cite{Mc1; Theorem 2} for the essential part of the argument; note that
Lemma 4 of \cite{MT} is incorrect if $U\cap\partial M\ne\emptyset$, for in
this case their $(\partial Q^3)-(\partial M^3)$ must have nonempty boundary
if $V$ is sufficiently small).
McMillan's proof is also based on Haken's Finiteness Theorem; in fact, our
proof of Theorem 2.8 can be used to simplify his argument by omitting the
reference to Waldhausen's result from the proof of \cite{Mc1; Theorem 2}.

\remark{Remarks. (i)}
Corollary 2.12 fails if compactness of $M$ is weakened to the assumption that
$M\but\partial M$ is the interior of a compact $3$-manifold and $\partial M$
is the interior of a compact $2$-manifold.
Figure 3 depicts a wedge of a PL knot $V$ (dashed) and a disguised version of
the Fox--Artin wild arc $Z\i S^3$ from Example 1.3 of \cite{FA}; the only
points of wildness of $Z$ are its endpoints $\{p,q\}$.
If $J$ is a closed subarc of $Z\but\{p,q\}$, its exterior in $Z$ consists of
two wild arcs $Y_p$ and $Y_q$, and it follows from \cite{FA; Example 1.2}
that $S^3\but (Y_p\cup Y_q)$ is homeomorphic to $S^3\but\{p,q\}$.
Let $R_p\sqcup R_q$ be a regular neighborhood of $Y_p\cup Y_q\but\{p,q\}$ in
$S^3\but\{p,q\}$ rel $\partial$.
The exterior $M$ of $R_p\cup R_q$ in $S^3\but\{p,q\}$ is a partial
compactification of $S^2\x\R$ with $\partial M\cong\R^2\sqcup\R^2$.

The compact subset $J\cup V$ of $M$ is contained in the contractible open
subset $U:=M\but X$, where $X$ is any PL arc in $S^3$, connecting $p$ and $q$
and disjoint from $J\cup V$ and from $R_p\cup R_q$.
(Indeed, arguments similar to those in \cite{FA; Example 2} show that
$U\but\partial U$ is homeomorphic to $S^3\but X$.)
On the other hand, suppose that $J\cup V$ is contained in a PL ball
$B^3$ in $M$.
Then $B^3$ meets $Z$ precisely in $J$, so $\partial B^3\but Z$ is a twice
punctured $2$-sphere.
Let $\gamma_0$ be a generator of $\pi_1(\partial B^3\but Z)\simeq\Z$.
Since $\pi_1(S^3\but Z)=1$ by \cite{FA}, but no power of the meridian
$\gamma_0$ is trivial in the tame knot group $\pi_1(B^3\but Z)$ (even modulo
commutators), by Seifert--van Kampen $\gamma_0$ has to be trivial in
$\pi_1(\Cl{S^3\but B^3}\but Z)$.
Let $W$ be any pushoff of $V$, disjoint from $Z$, such as the one in
\cite{FA; Fig\. 8}.
Then without loss of generality $W\i B^3$, so $\gamma_0$ is trivial in
$\pi_1(S^3\but (W\cup Z))$, contradicting \cite{FA}.
\endremark

\fig 3

\remark{(ii)} The set $J\cup V$ in Fig\. 3 is the image of a PL map
$f\:(I,\partial I)\to (M,\partial M)$ with one double point.
This $f$ is not a strong $\omega_1$-quasi-embedding (i.e\. there exists no PL
$3$-ball $B^3$ in $M$ such that $f^{-1}(B^3)$ is an arc) by the preceding
discussion.
However, $f$ is a strong $\omega$-quasi-embedding.
Indeed, the manifold $U$ from (i) is a partial compactification of $\R^3$
with $\partial U\cong\R^2\sqcup\R^2$.
For $1\le j<\omega$ let $P_j$ be any PL ball in $U\but\partial U$
containing the compact set $P_{j-1}\cup f(J_{j-1})$, and let $J_j$ be
the smallest subarc of $I$ containing $f^{-1}(P_j)$.
The assertion now follows from Proposition 2.1(a).
\endremark

\proclaim{Theorem 2.13} There exists an open $3$-manifold containing two
PL knots that are strongly $\omega$-quasi-isotopic but not PL isotopic.
\endproclaim

\demo{Proof} We use the notation of Remark (i) above.
Let $h\:\partial R_p\to\partial R_q$ be a homeomorphism between the boundary
planes of $M$ that identifies the endpoints of $J$.
The proper arcs $J$ and $J'$ in $M$, shown in Figures 3 and 4 respectively
(where $\partial J'=\partial J$) descend to knots $K$ and $K'$ in the open
$3$-manifold $M_h:=M/h$.
Since $M$ is simply-connected, there exists a generic PL homotopy
$H_t\:I\to M$ such that $H_0(I)=J$, $H_1(I)=J'$ and
$H_t(\partial I)=H_0(\partial I)$ for each $t\in I$.
By the argument of Remark (ii) above, $H_t$ descends to a strong
$\omega$-quasi-isotopy between $K$ and $K'$.

\fig 4

If $K_1$ is a PL knot in $M_h$, let $\Gamma(K_1)$ denote the kernel of
the inclusion induced homomorphism $\pi_1(M_h\but K_1)\to\pi_1(M_h)$.
If $K_2$ is obtained from $K_1$ by insertion of a local knot, the argument
in Remark (i) above shows that $\Gamma(K_1)=1$ iff $\Gamma(K_2)=1$.
Now $\Gamma(K)=1$ since a meridian of $J$ is obviously null-homotopic in $M$.
So to prove that $K$ and $K'$ are not PL isotopic it suffices to show that
$\Gamma(K')\ne 1$.
The arc $J'$ was chosen so that the mirror-symmetric halves of
the representation in \cite{FA; Example 1.3} extend to a nontrivial
representation $\rho\:\pi_1(M\but J')\to A_5$ in the alternating group;
the images of the additional Wirtinger generators are indicated in Fig\. 4.
Moreover, $\rho$ factors as $\rho\:\pi_1(M\but J')\to\pi_1(M_h\but K')\to A_5$,
where the first homomorphism is induced by the quotient map and the second
trivializes an additional generator of the HNN-extension $\pi_1(M_h\but K')$
of $\pi_1(M\but J')$.
Since $M\but J'$ is acyclic, $\pi_1(M_h\but K')$ is non-abelian.
But $\pi_1(M_h)=\Z$, so $\Gamma(K')\ne 1$. \qed
\enddemo

\head 3. $k$-cobordism and Cochran's invariants \endhead

\subhead \oldnos3.\oldnos1. Preliminaries \endsubhead
The results of this section are based on the following

\proclaim{Proposition 3.1} Suppose that links $L$ and $L'$ differ by a single
crossing change on the $i^{\text{th}}$ component so that the intermediate
singular link $f$ is a $k$-quasi-embedding.
Let $\ell$ denote the lobe $J_0$ of $f$, let $\hat\mu\in\pi(f)$ be a meridian
of $\ell$ which has linking number $+1$ with $\ell$, and let
$\hat\lambda\in\pi(f)$ denote the corresponding longitude of $\ell$
{\rm (see definition in \cite{MR})}.
Let $\mu,\lambda$ denote the images of $\hat\mu,\hat\lambda$ in $\pi(L)$, and
set $\tau=\mu^{-l}\lambda$, where $l$ is the linking number of the lobes.
Then $$\text{(a) \cite{MR; Lemma 3.1} } \lambda\in\underbrace{
\left<\mu\right>^{\left<\mu\right>^{\cdot^{\cdot^{\cdot^{\left<\mu\right>^{
\pi(L)}}}}}}\!\!;\!\!\!\!\!\!\!\!\!\!\!}_{k\text{ of }\,\left<\mu\right>
\text{'s}}\qquad\qquad \text{(b) }\tau\in[\dots[[\pi(L),\underbrace{
\left<\mu\right>],\left<\mu\right>]\dots,\left<\mu\right>}
_{k\text{ of }\,\left<\mu\right>\text{'s}}].$$
\endproclaim

Here $\pi(L)$ denotes the fundamental group $\pi_1(S^3\but L(mS^1))$, and
for a subgroup $H$ of a group $G$, the normal closure
$\left<g^{-1}hg\mid h\in H,g\in G\right>$ is denoted by $H^G$.
Part (b) can be proved analogously to (a), or deduced from it as follows.

\demo{Proof of (b)} Let us denote the subgroups in (a) and (b) by $A_k$ and
$B_k$, respectively.
Suppose that $\lambda=\mu^{n_1g_1}\dots\mu_r^{n_rg_r}$, where
each $n_i\in\Z$ and each $g_i\in A_{k-1}$.
Then $n_1+\dots+n_r=l$, so using the identity $ab=b[b,a^{-1}]a$ we can write
$$\mu^{-l}\lambda=\mu^{-n_1}\mu^{n_1g_1}h_1\mu^{-n_2}\mu^{n_2g_2}\dots
h_{r-1}\mu^{-n_r}\mu^{n_rg_r},$$ where
$h_i=[\mu^{n_ig_i},\mu^{l-(n_1+\dots+n_i)}]$.
But each $\mu^{-n_i}\mu^{n_ig_i}=[\mu^{n_i},g_i]\in B_k$ and each
$h_k^{-1}\in B_{k+1}\i B_k$ by the following lemma. \qed
\enddemo

\proclaim{Lemma 3.2} \cite{Ch; Lemma 2} Let $G$ be a group, $g\in G$, and
$n\in\N$.
Then $$[\left<g\right>,\underbrace{\left<g\right>^{\left<g\right>^{\cdot^{
\cdot^{\cdot^{\left<g\right>^G}}}}}\!\!\!\!\!}_{n\text{ of }\,\left<g\right>
\text{'s}}\ \,]=[\dots[[G,\underbrace{\left<g\right>],\left<g\right>]\dots,
\left<g\right>}_{n+1\text{ of }\,\left<g\right>\text{'s}}].$$
\endproclaim

Since we could not find an English translation of \cite{Ch}, we provide a short
proof for convenience of the reader.

\demo{Proof} By induction on $n$.
Since $[b,a]=[a,b]^{-1}$, $$\underbrace{
[[\dots[}_{n+1}G,\left<g\right>]\dots,\left<g\right>],\left<g\right>]=
[\left<g\right>,[\left<g\right>,\dots[\left<g\right>,G\underbrace{]\dots]]}
_{n+1}.$$
So it suffices to prove that $[\left<g\right>,[\left<g\right>,A_{n-1}]]=
[\left<g\right>,A_n]$, where $A_n$ denotes the expression in the left hand
side of the statement of the lemma, after the comma.
Indeed, pick any $h_1,\dots,h_r\in A_{n-1}$ and
$m,r,m_1,\dots,m_r,n_1,\dots,n_r\in\Z$, then
$$[g^m,[g^{m_1},h_1]^{n_1}\dots [g^{m_r},h_r]^{n_r}]=
[g^m,(g^{-m_1}g^{m_1h_1})^{n_1}
\dots (g^{-m_r}g^{m_rh_r})^{n_r}].$$
For some $s,m'_1,\dots,m'_s,m''_0,\dots,m''_s\in\Z$ and some
$h'_1,\dots,h'_s\in A_{n-1}$, the latter expression can be rewritten as
$$\multline[g^m,g^{m''_0}g^{m'_1h'_1}\dots
g^{m''_{s-1}}g^{m'_sh'_s}g^{m''_s}]=
[g^m,g^{m''_0}g^{m'_1h'_1}\dots
g^{m''_{s-1}+m''_s}g^{m'_sh'_sg^{m''_s}}]\\
=\dots=[g^m,g^{m''_0+\dots+m''_s}g^{m'_1h'_1 g^{m''_1+\dots+m''_s}}
\dots g^{m'_sh'_sg^{m''_s}}]=
[g^m,g^{m'_1h''_1}\dots g^{m'_sh''_s}]
\endmultline$$
for some new $h''_1,\dots,h''_s\in A_{n-1}$.
Clearly, this procedure is reversible. \qed
\enddemo

\subhead \oldnos3.\oldnos2. $k$-cobordism \endsubhead
We can now establish a relation with $k$-cobordism of Cochran \cite{Co3} and
Orr \cite{Orr1}.
We recall that the {\it lower central series} of a group $G$ is defined
inductively by $\gamma_1 G=G$ and $\gamma_{k+1}G=[\gamma_k G,G]$.
If $V$ is a properly embedded compact orientable surface in $S^3\x I$,
the {\it unlinked pushoff} of $V$ is the unique section $v$ of the spherical
normal bundle of $V$ such that $(v|_C)_*\:H_1(C)\to H_1(S^3\x I\but C)$
is zero for each component $C$ of $V$.

Two links $L_0,L_1\:mS^1\emb S^3$ are called {\it $k$-cobordant} if they can
be joined by $m$ disjointly embedded compact oriented surfaces
$V=V_1\sqcup\dots\sqcup V_m\i S^3\x I$ (meaning that
$V_j\cap S^3\x\{i\}=L_i(S^1_j)$ for $i=0,1$ and each $j$) so that if
$v\:V\emb S^3\x I\but V$ denotes the unlinked pushoff of $V$ then
the image of $v_*\:\pi_1(V)\to\pi_1(S^3\x I\but V)$ lies in the subgroup
generated by $v_*(\pi_1(\partial V))$ and $\gamma_k\pi_1(S^3\x I\but V)$.

\proclaim{Theorem 3.3} $k$-quasi-isotopy implies $(k+1)$-cobordism.
\endproclaim

\demo{Proof}
Let $h_t\:mS^1\to S^3$ be a $k$-quasi-isotopy, viewed also as
$H\:mS^1\x I\to S^3\x I$.
The (combinatorial) link of each double point $p$ of $H$ in the pair
$(S^3\x I,H(mS^1\x I))$ is a copy of the Hopf link in $S^3$.
Let us replace the star of $p$ in $H(mS^1\x I)$, which is the image of two
disks $D^2\x S^0$ in the same component of $mS^1\x I$, by an embedded
twisted annulus $A_p\simeq S^1\x D^1$ cobounded in $S^3$ by the components
of the Hopf link.
This converts $H(mS^1\x I)$ to an embedded surface $V\i S^3\x I$ with $m$
components and genus equal to the number of double points of $H$.
We may assume that the cylinder $A_p\simeq S^1\x I$ corresponding to
a double point $p=(t_p,x_p)\in S^3\x I$ meets the singular level
$S^3\x\{t_p\}$ in two generators $\{a,b\}\x I$, moreover $\{a\}\x I$
has both ends on the lobe $\ell_p=J_0$ of $h_{t_p}$.
Now $\pi_1(V)$ is generated by $\pi_1(\partial V)$ and the homotopy classes
of the loops $\tilde\ell_p:=(\ell_p\but H(D^2\x S^0))\cup\{a\}\x I$ and
$l_p:=S^1\x\{\frac12\}\i A_p$, connected to the basepoint of $V$ by some
paths.

Since $l_p\i A_p\i S^3$ bounds an embedded disk in $p*S^3$ whose interior
is disjoint from $V$, the pushoff $v(l_p)$ is also contractible in
the complement to $V$.
By Proposition 3.1(b), the free homotopy class of $v(\tilde\ell_p)$ lies in
$\gamma_{k+1}\pi(h_{t_p-\eps})$ and therefore in
$\gamma_{k+1}\pi_1(S^3\x I\but V)$. \qed
\enddemo

\remark{Remarks. (i)}
When $k=1$, the above proof works for weak $k$-quasi-isotopy in place of
$k$-quasi-isotopy, since the loop $v(\tilde\ell_p)$ is null-homologous both
in the complement to $V_i$, where $p\in H(S^1_i\x I)$ (by the definition of
the pushoff $v$) and in the complement to $\bigcup_{j\ne i}V_i$
(by the definition of weak $1$-quasi-isotopy).
In the case of two-component links with vanishing linking number this is not
surprising, since the Sato--Levine invariant $\beta=\bar\mu(1122)$,
which is well-defined up to weak $1$-quasi-isotopy \cite{MR; \S2}, is
a complete invariant of $2$-cobordism in this case \cite{Sa}.
\endremark

\remark{(ii)} By Dwyer's Theorem (see \cite{Orr1; Theorem 5}),
$\pi(L)/\gamma_{k+2}$ is invariant under $(k+1)$-cobordism and hence under
$k$-quasi-isotopy.
But an easier argument shows that even the $(k+3)^{\text{rd}}$ lower
central series quotient is invariant under $k$-quasi-isotopy.
Indeed, the quotient of $\pi(L)$ over the normal subgroup
$$\mu_k\pi(L)=\left<[m,m^g]\mid m\text{ is a meridian of }L,\,g\in
\underbrace{\left<m\right>^{\left<m\right>^{\cdot^{\cdot^{\cdot^{\left<m
\right>^{\pi(L)}}}}}}\!\!\!\!\!\!\!\!\!\!\!}_{k\text{ of }\,\left<m\right>
\text{'s}}\ \ \quad\right>$$
is invariant under $k$-quasi-isotopy by \cite{MR; Theorem 3.2}, and
$\mu_k\pi(L)\i\gamma_{k+3}\pi(L)$ by Lemma 3.2.
\endremark
\medskip

\subhead \oldnos3.\oldnos3. $\bar\mu$-invariants \endsubhead
Applying the result of \cite{Lin}, we get from Theorem 3.3 the following
statement with $2k+2$ in place of $2k+3$.

\proclaim{Theorem 3.4} Milnor's $\bar\mu$-invariants \cite{Mi2} of length
$\le 2k+3$ are invariant under $k$-quasi-isotopy.
\endproclaim

To take care of the remaining case (length equals $2k+3$), we give a direct
proof of Theorem 3.4, which is close to the above partial proof but avoids
the reference to Dwyer's Theorem in \cite{Lin}.

\demo{Proof} Consider a $\bar\mu$-invariant $\bar\mu_I$ where
the multi-index $I$ has length $\le 2k+3$.
If $L$ is a link, let $D_I(L)$ denote the link obtained from $L$ by replacing
the $i^{\text{th}}$ component by $n_i$ parallel copies, labelled
$i_1,\dots,i_{n_i}$, where $n_i$ is the number of occurrences of $i$ in $I$.
Then $\bar\mu_I(L)=\bar\mu_J(D_I(L))$, where $J$ is obtained from $I$ by
replacing the $n_i$ occurrences of each index $i$ by single occurrences of
$i_1,\dots,i_{n_i}$ in some order \cite{Mi2}.
Since $\bar\mu_J$ is an invariant of link homotopy \cite{Mi1}, \cite{Mi2}, it
suffices to show that if $I$ has length $\le 2k+3$ and $L_0$, $L_1$ are
$k$-quasi-isotopic, then $D_I(L_0)$, $D_I(L_1)$ are link homotopic.

Given a $k$-quasi-isotopy $h_t$ between $L_0$, $L_1$, we will again convert
it into a $(k+1)$-cobordism $V$ between $L_0$, $L_1$.
(Although this cobordism will be isotopic to the one constructed in the proof
of Theorem 3.3, here we are interested in attaching the handles of $V$
instantaneously.)
To this end, we emulate each crossing $L_-\mapsto L_+$ in $h_t$ by taking
a connected sum of $L_-$ with the boundary of a punctured torus $T$ in
the complement of $L_-$, as shown in Fig\. 5 (the ribbon, forming a half of
the punctured torus, is twisted $l$ times around the right lobe in order to
cancel the linking number $l$ of the two lobes).
The resulting link $L_-\#\partial T$ is easily seen to ambient isotop onto
$L_+$; meanwhile the natural generators of $\pi_1(T)$ represent the conjugate
classes of $\tau$ and $\mu^{-1}\mu^\tau=[\mu^{-1},\tau^{-1}]$, where $\mu$
is the meridian and $\tau$ the twisted longitude of the right lobe.
By Proposition 3.1(b), $\tau$ lies in $\gamma_{k+1}\pi(L_-)$, thus the
resulting surface $V$ is indeed a $(k+1)$-cobordism.

Replacing each component $V_i$ of $V$ by $n_i$ pairwise unlinked pushoffs,
we get a $(k+1)$-cobordism between $D_I(L_0)$ and $D_I(L_1)$.
The pushoffs of each punctured torus $T$ as above can be taken in the
same level $S^3\x\{t\}$ where these codimension one submanifolds become
naturally ordered.
Let us, however, shift these $n_i$ pushoffs in $S^3\x\{t\}$ vertically to
different levels $S^3\x\{t+\eps\},\dots,S^3\x\{t+n_i\eps\}$ in the order just
specified.
Thus the $(k+1)$-cobordism splits into a sequence of isotopies and
iterated additions of boundaries of punctured tori:
$$D_I(L_-)\ \ \mapsto\ \ D_I(L_-)\underset\ i_1\!\to\#\partial T_1\ \ \mapsto
\ \ \dots\ \ \mapsto\ \ D_I(L_-)\underset\ i_1\!\to\#\partial T_1\underset
\ i_2\!\to\#\cdots\underset\ \,i_{n_i}\!\!\!\to\#\partial T_{n_i}\!=\!
D_I(L_+).$$
(The subscript of $\#$ indicates the component being amended.)
It remains to show that every two consecutive links
$L(j):=D_I(L_-)\#_{i_1}\dots\#_{i_j}T_{i_j}$ and $L(j+1)$ in such a string
are link homotopic.
Now $L(j)$ and $L(j+1)$ only differ in their $(i_{j+1})^{\text{th}}$
components $K(j)$ and $K(j+1)=K(j)\#\partial T_{i_{j+1}}$, and by \cite{Mi1}
it suffices to show that these represent the same conjugate class of Milnor's
group $\pi(L'_j)/\mu_0\pi(L'_j)$ of $L'_j:=L(j)\but K(j)=L(j+1)\but K(j+1)$,
in other words, that the conjugate class of $\partial T_{i_{j+1}}$ is
contained in the normal subgroup $\mu_0\pi(L'_j)$.

\fig 5

In the case $j=0$ we can simply quote the above observation that
the generators $\sigma,\tau$ of $\pi_1(T)$ lie in $\gamma_{k+1}\pi(L_-)$ and
$\gamma_{k+2}\pi(L_-)$, respectively.
Hence the class of $\partial T_{i_1}$ lies in
$\gamma_{2k+3}\pi(L'_0)\i\gamma_m\pi(L'_0)\i\mu_0\pi(L'_0)$, where the latter
containment holds by \cite{Mi1} since $L'_0=L_-\but K(0)$ has $m-1$
components.
On the other hand, in the case $j=n_i-1$ it suffices to notice that one
generator $\sigma$ of $\pi_1(T)$ bounds an embedded disk in the complement
to $L_+$; hence $\partial T_{i_{n_i}}$ is null-homotopic in the complement
to $L'_{n_i-1}=L_+\but K(n_i)$.
Finally, when $0<j<n_i-1$, it takes just a little more patience to check
that the generators $\tau_{i_{j+1}}$, $\sigma_{i_{j+1}}$ of $T_{i_{j+1}}$
are homotopic in the complement to $L(j)$ respectively to $\tau$ and
$(\mu_{i_{j+1}}\dots\mu_{i_{n_i}})^{-1}(\mu_{i_{j+1}}\dots\mu_{i_{n_i}})^\tau=
[(\mu_{i_{j+1}}\dots\mu_{i_{n_i}})^{-1},\tau^{-1}]$, where $\tau$ denotes
the longitude of the right lobe and $\mu_{i_l}$ the meridian of
the $(i_l)^{\text{th}}$ component of $L(j)$.
(To see this, one may assume that $n_i=2$, since all non-amended components
can be grouped together and all amended components can also be grouped
together, using that the punctured tori $T_{i_j}$ are ordered in agreement
with position of their projections to $S^3\x\{t\}$.)
Thus $\tau_{i_{j+1}}\in\gamma_{k+1}\pi(L(j))$ and
$\sigma_{i_{j+1}}\in\gamma_{k+2}\pi(L(j))$, whence
the class of $\partial T_{i_{j+1}}$ lies in
$\gamma_{2k+3}\pi(L'_j)\i\gamma_m\pi(L'_j)\i\mu_0\pi(L'_j)$.
\qed
\enddemo

\remark{Remark}
The restriction $2k+3$ is sharp, since Milnor's link $\M_{k+1}$ is strongly
$k$-quasi-isotopic to the unlink (see \S1) but has nontrivial
$\bar\mu(\underbrace{11\dots 11}_{2k+2}22)$ (cf\. \cite{Mi2}).
\endremark

\subhead \oldnos3.\oldnos4. Conway polynomial and Cochran's invariants
\endsubhead
We recall that the Conway polynomial of an $m$-component link is of the
form $z^{m-1}(c_0+c_1z^2+\dots+c_nz^{2n})$ (see, for instance, \cite{Co2}).

\proclaim{Corollary 3.5}
Let $c_k$ denote the coefficient of the Conway polynomial at $z^{m-1+2k}$.
The residue class of $c_k$ modulo $\gcd(c_0,\dots,c_{k-1})$ is invariant
under $k$-quasi-isotopy.

Set $\lambda=\lceil\frac{(l-1)(m-1)}2\rceil$.
The residue class of $c_{\lambda+k}$ modulo $\gcd$ of
$c_\lambda,\dots,c_{\lambda+k-1}$ and all $\bar\mu$-invariants of length
$\le l$ is invariant under
$(\lfloor\frac\lambda{m-1}\rfloor+k)$-quasi-isotopy.
\endproclaim

Here $\lceil x\rceil=n$ if $x\in [n-\frac12,n+\frac12)$, and
$\lfloor x\rfloor=n$ if $x\in (n-\frac12,n+\frac12]$ for $n\in\Z$.
The second part contains the first as well as Cochran's observation that
$c_1$ is a homotopy invariant if $m\ge 3$ and the linking numbers vanish
\cite{Co2; Lemma 5.2}.

\demo{Proof} Clearly, $k$-quasi-isotopic links are closures of
$k$-quasi-isotopic string links.
(String links and their closures are defined, for instance, in \cite{Le}, and
we assume their $k$-quasi-isotopy to be fixed on the boundary.)
It is entirely analogous to the proof of Theorem 3.4 to show that
$\mu$-invariants of string links of length $\le 2k+3$, as well as their
modifications in \cite{Le} are invariant under $k$-quasi-isotopy.
(The novelty in the definition of $\mu$-invariants in \cite{Le} is that
instead of the longitude $\lambda_i$ one uses the product
$\lambda_i\mu_i^{-l_i}$, where $\mu_i$ is the corresponding meridian, and
$l_i$ is the sum of the linking numbers of the $i$-th component with
the other components; in particular, this leads to $\mu(i,i)=-l_i$.)
Let $a_i(S)$ denote the coefficient at $z^i$ of the power series
$\Gamma_S(z)$, defined for a string link $S$ in the statement of Theorem 1
of \cite{Le}.
(Note that the determinant in the definition of $\Gamma_S$ refers to any
diagonal $(m-1)\x (m-1)$ minor, rather than to the full matrix.)
The definition of $\Gamma_S$ shows that $a_{m-1+2k}(S)$ is determined by
the $\mu$-invariants of $S$ of length $\le 2k+2$, and hence is invariant
under $k$-quasi-isotopy by Theorem 3.4.
On the other hand, \cite{Le; Theorem 1} (see also \cite{KLW}, \cite{TY},
\cite{MV} and \cite{Tr; Theorem 8.2 and the third line on p\. 254}) implies
that the closure $L$ of $S$ satisfies $c_k(L)\equiv a_{m-1+2k}(S)$ modulo
the greatest common divisor of $a_i(S)$ with $i<m-1+2k$.
Since the $c_k$'s are the only possibly nonzero coefficients of the Conway
polynomial, the same congruence holds modulo $\gcd(c_0,...,c_{k-1})$.

The second part follows by the same argument, taking into account that,
by the definition of $\Gamma_S$, the residue class of $a_{l(m-1)+n}$ modulo
$\gcd$ of all $\mu$-invariants of $S$ of length $\le l$ (equivalently, of all
$\bar\mu$-invariants of $L$ of length $\le l$) is determined by
the $\mu$-invariants of $S$ of length $\le l+n+1$. \qed
\enddemo

From Theorem 3.3 and \cite{Co4; 2.1} we immediately obtain

\proclaim{Corollary 3.6} Cochran's derived invariants $\beta^i$, $i\le k$,
of each two-component sublink with vanishing linking number are invariant
under $k$-quasi-isotopy.
\endproclaim

Using Cochran's original geometric definition of his invariants \cite{Co1}
(see also \cite{TY}) it is particularly easy to see that $\beta^k(\M_k)=1$
(whereas $\beta^k($unlink$)=0$).
Thus $\M_k$ is not $k$-quasi-isotopic to the unlink.
On the other hand, the statements of Corollary 3.6 and Theorem 3.3 are sharp,
that is, $\beta^{k+1}$ is not invariant under $k$-quasi-isotopy, which
therefore does not imply $(k+2)$-cobordism.
In fact, it is known \cite{Co3}, \cite{St} (see also \cite{Orr2}) that
$\beta^i$ is an integral lifting of
$\pm\bar\mu(\underbrace{11\dots 11}_{2i\text{ times}}22)$.

Cochran showed that the power series $\sum_{k=1}^\infty\beta^kz^k$ is
rational, and is equivalent to the rational function $\eta^+_L(t)$ of
Kojima--Yamasaki \cite{KY} (see definition in {\it ``$n$-Quasi-isotopy I''})
by the change of variable $z=2-t-t^{-1}$ \cite{Co1; ~\S7} (see also
\cite{TY}).
Kojima and Yamasaki were able to prove that the $\eta$-function is invariant
under topological I-equivalence (i.e\. the non-locally-flat version of
concordance) \cite{KY}.
However, the following comments are found in \cite{KY; Introduction}:
``In the study of the $\eta$-function, we became aware of impossibility to
define it for wild links.
The reason of this is essentially due to the fact that the knot module of
some wild knot is not $\Z[t^{\pm 1}]$-torsion.''
Indeed, in the definition of $\eta^+_L$, where $L=K_+\cup K_-$, one has to
divide by a nonzero Laurent polynomial annihilating the element $[\tilde K_+]$
of the knot module of $K_-$.
It is shown in \cite{KY} that no such Laurent polynomial exists for a twisted
version $\widetilde\M_\infty$ of Milnor's wild link $\M_\infty$ (which can be
approximated by a twisted version $\widetilde\M_n$ of the links $\M_n$, with
one half-twist along each disk as in Fig.\ ~1).

On the other hand, Corollary 3.6 and \cite{MR; Corollary 1.4(a)} imply

\proclaim{Corollary 3.7} Given any $i\in\N$ and any topological link
$\frak L$, there exists an $\eps>0$ such that all PL links, $C^0$
$\eps$-close to $\frak L$, have the same $\beta^i$.
\endproclaim

This has already been known, by Milnor's work on $\bar\mu$-invariants
\cite{Mi2}, for the residue class of $\beta^k$ modulo
$\gcd(\beta^1,\dots,\beta^{k-1})$.
Corollary 3.7 yields a natural extension of each $\beta^i$ to all
topological links, which by \cite{MR; Corollary 1.4(b)} or by compactness
of $I$ is invariant under topological isotopy, i.e\. homotopy through
embeddings.
In particular, since $\widetilde\M_\infty$ is isotopic to the unlink (see
\cite{Mi2}), Cochran's reparametrization $\sum\beta^iz^i$ of $\eta^+_L$ of
an arbitrary PL $\eps$-approximation of $\widetilde\M_\infty$ converges to
zero (as a formal power series) as $\eps\to 0$.

\head Acknowledgements \endhead

The first author would like to thank the University of Florida math faculty for
their hospitality and to P. Akhmetiev, T. Cochran, J. Male\v{s}i\v{c},
R. Mikhailov and R. Sadykov for useful discussions.
We are grateful to the referee for pointing out a number of inaccuracies in
the first version.
The authors are partially supported by the Ministry of Science and Technology
of the Republic of Slovenia grant No\. J1-0885-0101-98 and the Russian
Foundation for Basic Research grants 99-01-00009, 02-01-00014 and 05-01-00993.

\enddocument
\end